# Modeling Microstrip Antenna

Luis Alberto Rabanal Ramirez[1], *Cláudio Márcio de Freitas Silva*[2]
[1] Universidade Estadual do Norte Fluminense UENF, Brazil (prof.ramirez.luis@gmail.com).
[2] Centro Universitário de Volta Redonda UNIFOA, Brazil. (Claudio.freitas@foa.org.br).

*Abstract*—In this work, a rectangular microstrip antenna with inset is designed, simulated and optimized. In the optimization process the patch is deformed, it new antenna present a amorphous patch. The optimization process was conducted with Genetic Algorithm (GA), $S_{11}$ parameters was obtained with full wave Finite-Differences Time-Domain (FDTD-3D), and the initial configuration (design) was obtained with line transmission and cavite method. All methods implemented in-house software. The antenna was designed to operate in Ku band with the center frequency at 16GHz, the frequency band (15.25–17.25GHz). Antenna with return loss minor that -27 dB and bandwidth between 3-4 Ghz, is obtained. The effectiveness of the proposed designs is confirmed through proper simulation results.

## I. INTRODUCTION

The rapid advancement in the satellite communication field in the past few decades has led to the development of small, low profile and efficient antenna design. Antenna is an important structure in any satellite communication system, service as satellite internet access, spacecraft telemetry, command, and tracking comunications[1], they has boosted research in this field, and good antenna design definitely improves the overall performance of the system.

This work focuses in microstrip antennas for satellite communication, because it present several advantages compared to another antennas, by sample, large single-dish[1,2]. Various schemes are being used to minimize the drawbacks of this type of antennas, but it is a difficult task. In this work, the antenna was design to operate in Ku band with the center frequency at 16GHz, the frequency band (15.25–17.25GHz). This work used the genetic algorithms to overcome these disadvantages, based on references [3-5]. Specifically the numerical method FDTD-3D was combined with the genetic algorithm "GA" for optimization the bandwidth and return loss, The initial geometric dimensions was obtained used the ressonante cavity / line transmition methods [6], All methods implemented in-house software.

The FDTD-3D details relating to modeling the antenna was described in section 2. Section 3 presents aspect relationated with Genetic Algorithms used in this work. The section 4 shows the present numerical results. Finally, section 5 presents the conclusions and achieve some goals. The effectiveness of the proposed designs is confirmed through proper simulation results.

## II. FDTD-3D METHOD

The FDTD-3D method with the "Uniaxial Perfect Matched Layer" (UPML) [7] was used.

### A. The source

The Morlet wavelet function was used. The excitement given by the equation (1)

$$E_z(t) = E_0 \exp[-(2\pi f_b(t-t_0))^2]\cos[2\pi f_c(t-t_0)] \quad (1)$$

where $E_0$ = 22 V/m is the signal amplitude, $t_0$ = 5.02 ns is the time at which the pulse reaches its maximum value, fc = 16 GHz is the central frequency, fb = 2000 MHz and $0 \leq t \leq 8000\Delta t$ represents the propagation time signal with Δt = 0.0342 ns representing the value of the time step.

Figure 1 shows the time domain signal and its frequency domain equivalent.

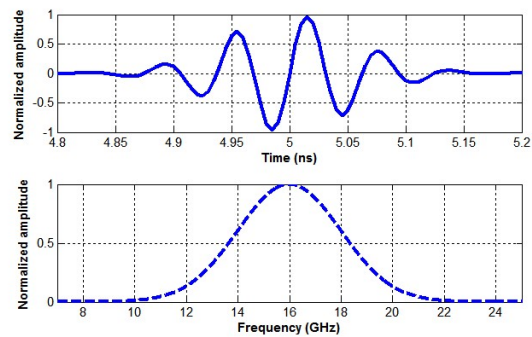

Fig. 1. wavelet function for excitation: (a) time domain and (b) frequency spectrum.

### B. Scattering parameter

The traditional procedure for obtaining the corresponding frequency domain, aiming to calculate the return loss $S_{11}$ of the microstrip patch is know the incident and reflected fields. The total field in FDTD computational domain is the sum of the incident and reflected fields. To achieve the incident field suppress from the simulation the patch that cause reflections and is calculate the fields in line of infinite length. This incident field is subtracted from the field generated by the FDTD, resulting in the reflected field. $S_{11}$ obtains the parameter (ω) by applying the Fourier transform of the reflected field. This procedure entails running the program twice for the parameters $S_{11}$, resulting in a lengthy process, new schemes have emerged as those that solve linear equations [8], others that improve the process of transformation of fourier [9]. In this work when the evolutionary process begins (specifically the first time) using a matrix help to save the incident field values before interacting with the patch as soon as these data are used to calculate the field reflected in future generations.

## III. GENETIC ALGORITMH

The optimization process begins considering the rectangular microstrip antenna with inset geometrically scaled using the methods transmission line and resonant cavity[6] (Figure 2).

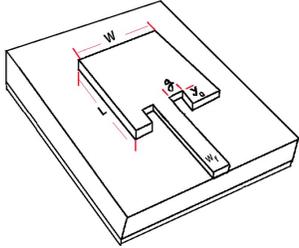

Fig. 2. Microstrip antenna with inset (unscaled)

The parametros values are: W=11.70mm, L=8.87mm, $W_f$=1.40mm, g=0.7mm, $Y_0$=1.33mm, material substrate value = 2.2, height substrate =0.76 mm, the ground plane and height patch and line were taken as one delta mesh in the z direction.

This patch was divided into a binary matrix of 17 x 17 sub-patches, as show in Figure 3. The sub-patch in this representation is an intrinsic representation of problem, this indicates presence (set as 1) or absence (set as 0) of copper sub-patches in candidate solution (chromosome).

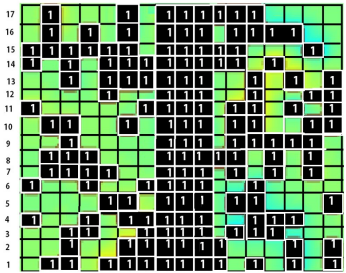

Fig. 3. pixelated amorphous patch antenna.

The population size was choose in 30, with relation to chromosome size, various sizes were tested, but a digital individual of 289 elements (connect rows of the matrix).

The fitness function is related to the bandwidth and the return loss, equation (2). The two optimization objectives functions were combined to form a single overall fitness function.

$$fit = \alpha\, fit\_BW + (1-\alpha)\, fit\_PR \qquad (2)$$

in (2) α is a weighting factor in the range [0, 1] this parameter allows us to choose the emphasis or preference of a parameter with respect to each other. The fitness function related to bandwidth is set with minimum bandwidth 2,0 GHz, and the related to return loss demand further loss of -25 dB.

After the initialization of population and defined the fitnees function, begin the cyclic process or evaluative process, cycle and generation are effectively interchangeable terms. In the cycle process was used roulette-wheel technique; shuffle crossover with swapping probability equal to 0.5; for the mutation probability was fixed in 0.001 and bit-flip[10].

## IV. NUMERICAL RESULT

An antenna rectangular microstrip with inset, to work in the Ku band, whose dimensions were obtained by the methods resonant cavity and transmission line, simulated with FDTD-3D, it has been optimized in order to provide higher bandwidth (compared to original), using genetic algorithms, one schemes for the representation of geometry were used, which generate one optimized configuration with reference the original.

Figure (4) the line feed was extended to the top edge of the patch, and the current will flow from the line to joined sub-patch (form a resonant structure). These sub-patch that not directly couple became a parasitic sub-patches. In regular antennas is know, that the present of these elements enhanced the wideband of the antenna, in Figure (4.b), the wideband is 3.0 GHz a return loss is -34.84 dB.

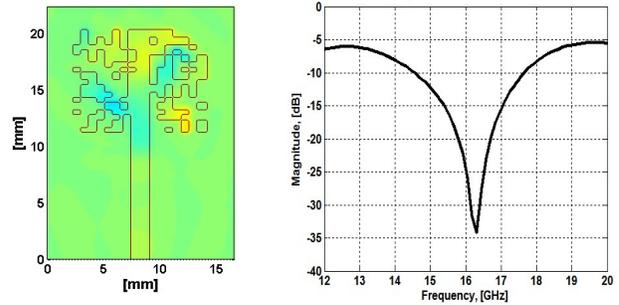

Fig. 4. Intermediate results of the optimization process. (a) amorphouse shape antenna, (b) return loss graph.

Figure (5) present more return loss if compared with Figure (4), but the frequency of resonance departs 16 GHz, chosen as the central frequency, minor band-width, aproximately 2 GHz.

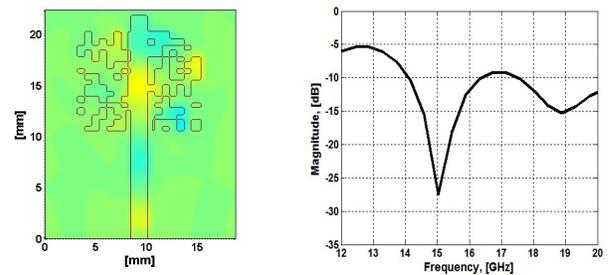

Fig. 5. Intermediate results of the optimization process. (a) amorphous shape antenna, (b) return loss graph.

these changes would be due to mutations, in particular the extension beyond the upper edge of the patch, feed line (in 2 subpatchs). it is also possible to observe in this figure that the return loss (in the range between 16.2 and 17.3 GHz) was less than -10 dB, the bandwidth for this configuration would be higher than 6 GHz.

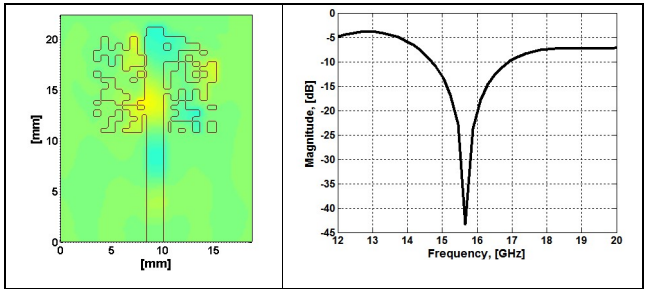

Fig. 6. final results of the optimization process. (a) amorphouse shape antenna with notch and the end of line fed , (b) correspondent return loss.

Figure (6) present the final result of the optimization process, in (6.a) the line feed was extended beyond to the top edge of the patch (in 1 subpatch), the wideband of the antenna, in Figure (6.b), the wideband is 2.2 GHz a return loss is -44.26 dB. Small changes could be inserted, like a manual mutation or fine adjustment of the mutation, based on techniques to increase bandwidth (amply documented).

## V. CONCLUSION

The aim of this work was to design a compact, efficent microstrip patch antenna for use in Ku band. As can be seen in the figures and (4-6), one obtains return values and loss of bandwidth considered good, if compared with the rectangular antenna with inset used as the original chromosome.

The conclusion general is that after the optimization process non restricted configurations are possible. Consequently, they can potentially discover beneficial structures that cannot be discovered by conventional design approaches.